\documentclass{article}
\usepackage{amsfonts}
\usepackage{amssymb}
\usepackage{graphicx}
\usepackage{amsmath}

\input{tcilatex}
\begin{document}

\begin{center}
\textbf{Veli\ B. Shakhmurov}\ 

\textbf{\ Completness of roots elementes of linear operators in Banach
spaces and application}

\ \ \ \ \ \ \ \ \ \ \ \ \ \ \ \ \ \ \ \ \ \ \ \ \ \ \ \ \ \ \ \ \ \ \ \ \ \
\ \ \ \ \ \ \ \ \ \ \ 

\bigskip Okan University, Department of Mechanical Engineering , Akfirat,
Tuzla 34959 Istanbul, Turkey, E-mail: veli.sahmurov@okan.edu.tr

\ \ \ \ \ \ \ \ \ \ \ \ \ \ \textbf{ABSTRACT}
\end{center}

In this paper the general spectral properties of linear operators in Banach
spaces are studied. We find sufficient conditions on structure of Banach
spaces and resolvent properties that guarantee completeness of roots
elements of Schatten class operators. This approach generalizes the well
known result for operators in Hilbert spaces. In application, the boundary
value problems for the abstract equation of second order with variable
coefficients are studied. The principal part of the appropriate differential
operator is not self-adjoint. The discreetness of spectrum and completeness
of root elements of this operator are obtained.

\textbf{Key Words}{: }Uniformly convex Banach spaces; Abstract functions;
Schatten class of operators; Completeness of root elements; Separable
boundary value problems; Differential-operator equations;

\textbf{AMC 2000: 47Axx, 47A10, 35Jxx, 35Pxx \qquad\ }

One of the fundamental results on spectral theory of operators is the
completeness of roots elements of Schatten class operators in Hilbert spaces:

\textbf{Theorem} $\left[ \text{8, Theorem XI. 9.29}\right] $. Assume:

(1) $H$ is a Hilbert space and $A$ is an operator in $C_{p}\left( H\right) $%
, for some $p\in \left( 1,\infty \right) ;$

(2) $\gamma _{1},\gamma _{2},...,\gamma _{s}$ is non overlapping,
differentiable arcs in the complex plane starting at the origin. Suppose
that each of the $s$ regions into which the plans is divided by these arcs
is contained in an angular sector of opening less then $\frac{\pi }{p}.$
Moreover, $m>0$ is an integer so that the resolvent of $A$ satisfies the
inequality $\left\Vert R\left( \lambda ,A\right) \right\Vert =O\left(
\left\vert \lambda \right\vert ^{-m}\right) $ \ as $\lambda \rightarrow 0$
along any of the arcs $\gamma _{i}$. Then the $spA$ contains the subspace $%
A^{m}H.$

The main aim of the present paper is the generalization of the above
important theorem $\left[ \text{8, Theorem XI. 9.29}\right] $ for Banach
spaces. The spectral properties of linear operators in Banach spaces is a
subject which is not much investigated. The related effort, indeed requires
new tools of modern analysis and operator theory. Nevertheless, the results
in this field have numerous applications in pure differential, pseudo
differential and functional-differential equations. For this reason, it was
very important to have general result about spectral properties of linear
operators in Banach spaces. The articles $\left[ 2\right] $, $\left[ 6\right]
$ and $\left[ 15\right] $ are devoted to this question in Banach spaces. In
this paper, we disclose different sufficient condition for completeness of
roots elements of linear operators. We consider the class of Banach spaces
which satisfy some given conditions, but by virtue of Remark1, our class of
operators are wider than the class of operators considered in $\left[ 2%
\right] $, $\left[ 6\right] $ and $\left[ 15\right] .$ Also, in $\left[ 6%
\right] $ the extra condition is assumed to be nonempty of spectrum of these
class of operators. Moreover, our method of proofs are different from proofs
provided in the cited references $.$

We find sufficient conditions on structure of Banach spaces which allow to
define the trace of operators and its properties. Also, we get Carleman
estimate of quasi nuclear operators (QNOs) and its specral properties. In
application we consider nonlocal boundary value problem (BVP)\ for the
second order differential-operator equation (DOE) with top variable
coefficients 
\begin{equation*}
Lu=a\left( x\right) u^{\left( 2\right) }\left( x\right) +B\left( x\right)
u^{\left( 1\right) }\left( x\right) +A\left( x\right) u\left( x\right)
=f\left( x\right) ,\text{ }x\in \left( 0,1\right) ,
\end{equation*}%
where $a_{k}$ are complex-valued functions, $A\left( x\right) $, $B\left(
x\right) $ are linear operators in a Banach space $E$ and $f$ is a $E$%
-valued function$.$ The principal part of the associate differential
operator is not self-adjoint. We prove that, the spectrum of the associated
differential operator is discrete and the system of roots elements are
complete in $E$-valued weighted $L_{p}$ spaces. Note that,
differential-operator equations (DOEs) have been studied extensively by many
researchers (see $\left[ \text{1, 3}\right] ,$ $\left[ \text{7, 9, 11, 13, 14%
}\right] ,$ $\left[ \text{16-26}\right] $ and the references therein).

We start by giving the notations and definitions to be used in this paper.

\ Let $\gamma =\gamma \left( x\right) $ be a positive measurable weighted
function on the region $\Omega \subset R^{n}$. Let $L_{p,\gamma }\left(
\Omega ;E\right) $ denote the space of all strongly measurable $E$-valued
functions that are defined on $\Omega $ with the norm

\begin{equation*}
\left\Vert f\right\Vert _{p,\gamma }=\left\Vert f\right\Vert _{L_{p,\gamma
}\left( \Omega ;E\right) }=\left( \int \left\Vert f\left( x\right)
\right\Vert _{E}^{p}\gamma \left( x\right) dx\right) ^{\frac{1}{p}},1\leq
p<\infty \ .
\end{equation*}

The weight $\gamma \left( x\right) $ we will consider satisfy an $A_{p}$
condition. i.e., $\gamma \left( x\right) \in A_{p},$ $p\in \left( 1,\infty
\right) $ if there is a positive constant $C$ such that

\begin{equation*}
\left( \frac{1}{\left\vert Q\right\vert }\dint\limits_{Q}\gamma (x)dx\right)
\left( \frac{1}{\left\vert Q\right\vert }\dint\limits_{Q}\gamma ^{-\dfrac{1}{%
p-1}}(x)dx\right) ^{p-1}\leq C
\end{equation*}%
for all balls $Q\subset R^{n}$.

For $\gamma \left( x\right) \equiv 1$ the space $L_{p,\gamma }\left( \Omega
;E\right) $ will be denoted by $L_{p}\left( \Omega ;E\right) .$ The Banach
space $E$ is said to be a $\zeta $-convex space (see e.g. $\left[ 4\right] $%
) if there exists a symmetric real-valued function $\zeta \left( u,v\right) $
on $E\times E$ which is convex with respect to each of the variables, and
satisfies the conditions 
\begin{equation*}
\zeta \left( 0,0\right) >0,\text{ }\zeta \left( u,v\right) \leq \left\Vert
u+v\right\Vert \text{ for }\left\Vert u\right\Vert =\left\Vert v\right\Vert
=1.
\end{equation*}%
\ The Banach space\ $E$ is called an $UMD$-space if\ the Hilbert operator $%
\left( Hf\right) \left( x\right) =\lim\limits_{\varepsilon \rightarrow
0}\int\limits_{\left\vert x-y\right\vert >\varepsilon }\frac{f\left(
y\right) }{x-y}dy$ \ is bounded in $L_{p}\left( -\infty ,\infty ,E\right) ,$ 
$p\in \left( 1,\infty \right) $ (see. e.g. $\left[ 4\right] $ ). $UMD$
spaces include e.g. $L_{p}$, $l_{p}$ spaces and Lorentz spaces $L_{pq},$ $p$%
, $q\in \left( 1,\infty \right) $. It is shown $\left[ 4\right] $ that the
Banach space $E$ is $UMD$ if only if this space is a $\zeta $-convex space.

Let $\mathbb{C}$ be the set of complex numbers and

\ 
\begin{equation*}
S_{\varphi }=\left\{ \lambda ;\text{ \ }\lambda \in \mathbb{C}\text{, }%
\left\vert \arg \lambda \right\vert \leq \varphi \right\} \cup \left\{
0\right\} ,0\leq \varphi <\pi .
\end{equation*}

Let $E_{1}$ and $E_{2}$ be two Banach spaces. $B\left( E_{1},E_{2}\right) $
denotes the space of bounded linear operators from $E_{1}$ to $E_{2}$. For $%
E_{1}=E_{2}=E$ it will be denoted by $B\left( E\right) .$

A linear operator\ $A$ is said to be positive in a Banach\ space $E$ with
bound $M>0$ if $D\left( A\right) $ is dense on $E$ and 
\begin{equation*}
\ \left\Vert \left( A+\lambda I\right) ^{-1}\right\Vert _{B\left( E\right)
}\leq M\left( 1+\left\vert \lambda \right\vert \right) ^{-1}
\end{equation*}%
with $\lambda \in S_{\varphi },\varphi \in \left( 0,\pi \right] $, $I$ is an
identity operator in $E.$ Sometimes instead of $A+\lambda I$\ will be
written $A+\lambda $ and it will be denoted by $A_{\lambda }$. Let $E\left(
A\right) $ denote $D\left( A\right) $ with the graphical norm.

A set $W\subset B\left( E_{1},E_{2}\right) $ is called $R$-bounded (see e.g, 
$\left[ \text{9}\right] $ ) if there is a constant $C>0$ such that for all $%
T_{1},T_{2},...,T_{m}\in W$ and $u_{1,}u_{2},...,u_{m}\in E_{1},m\in N$ 
\begin{equation*}
\int\limits_{0}^{1}\left\Vert \sum\limits_{j=1}^{m}r_{j}\left( y\right)
T_{j}u_{j}\right\Vert _{E_{2}}dy\leq C\int\limits_{0}^{1}\left\Vert
\sum\limits_{j=1}^{m}r_{j}\left( y\right) u_{j}\right\Vert _{E_{1}}dy,
\end{equation*}%
where $\left\{ r_{j}\right\} $ is a sequence of independent symmetric $%
\left\{ -1,1\right\} $-valued random variables on $\left[ 0,1\right] $.

The positive operator $A$ is said to be an $R$-positive in a Banach space $E$
if there exists $\varphi \in \left[ 0\right. ,\left. \pi \right) $ such that
the set $\left\{ A\left( A+\xi I\right) ^{-1}:\xi \in S_{\varphi }\right\} $
\ is $R$-bounded.

A linear operator\ $A=A(x),$ $x\in \left[ a,b\right] $ is said to be
uniformly positive in a Banach\ space $E$, if $D\left( A(x)\right) $ dense
in $E\ $and does not depend on $x$ and there is a constant $M>0$ such that 
\begin{equation*}
\ \left\Vert \left( A(x)+\lambda I\right) ^{-1}\right\Vert _{B\left(
E\right) }\leq M\left( 1+\left\vert \lambda \right\vert \right) ^{-1}\text{ }
\end{equation*}
for all $\lambda \in S_{\varphi }$, $x\in \left[ a,b\right] $ and some $%
\varphi \in \left[ 0\right. ,\left. \pi \right) .$

Let $E_{0}$ and $E$ be two Banach spaces and $E_{0}$ is continuously and
densely embeds into $E$.

Let\ $W_{p,\gamma }^{2}\left( 0,1;E_{0},E\right) $ denote a space of all
functions $u\in L_{p,\gamma }\left( 0,1;E_{0}\right) $ possess the
generalized derivatives $u^{\left( 2\right) }\in L_{p,\gamma }\left(
0,1;E\right) $ with the norm 
\begin{equation*}
\ \left\Vert u\right\Vert _{W_{p,\gamma }^{2}\left( 0,1;E_{0},E\right)
}=\left\Vert u\right\Vert _{L_{p,\gamma }\left( 0,1;E_{0}\right)
}+\left\Vert u^{\left( 2\right) }\right\Vert _{L_{p,\gamma }\left(
0,1;E\right) }<\infty .
\end{equation*}

Sp $A$ denote the closure of the linear span of the roots elements of the
operator $A.$

Let $E$ be a Banach space and $E^{\ast }$ denotes its dual. For $u\in E,$ $%
f\in E^{\ast }$ let $<u,f>$ denote the value of $f$\ \ for $u$, i.e. $%
<u,f>=f\left( u\right) $. Suppose $\left\{ e_{j},f_{j}\right\} ,$ $j=1,2,...$
is a biorthonormal basis systems in $E\times $ $E^{\ast }$, i.e. 
\begin{equation*}
\left\{ e_{j}\right\} \subset E,\left\{ f_{j}\right\} \subset E^{\ast },%
\text{ }<e_{j},f_{i}>=\delta _{ij},\text{ }\left\Vert e_{j}\right\Vert
_{E}=1,\text{ }\left\Vert f_{i}\right\Vert _{E^{\ast }}=1,\text{ }i,\text{ }%
j=1,2,....
\end{equation*}

For $u\in E,$ $f\in E^{\ast }$ let $\alpha _{j}=$ $<u,f_{j}>$ and $\beta
_{j}=$ $<e_{j},f>$ denote the Fourier coefficients of $u$ and $f$ with
respect to systems $\left\{ e_{j}\right\} \subset E$ and $\left\{
f_{j}\right\} \subset E^{\ast }$, respectively.

\textbf{Definition1. }A separable Banach space with base is said to be the
space satisfying the $B$-condition, if there are a positive constant $C$ and
a $p$ $\in \left( 1,\infty \right) $ such that 
\begin{equation*}
<u,f>=\sum\limits_{j=1}^{\infty }\alpha _{j}\beta _{j}\text{, }\left\Vert
u\right\Vert _{E}^{p}\leq C\sum\limits_{j=1}^{\infty }\left\vert \alpha
_{j}\right\vert ^{p}<\infty
\end{equation*}%
for all biorthonormal basis systems $\left\{ e_{j},f_{j}\right\} ,$ $%
j=1,2,...$ in $E\times $ $E^{\ast }.$

The Hilbert spaces satisfies this condition for $p=2.$ For examples $L_{p}$
and $l_{p}$ spaces, $p\in \left( 1,\infty \right) $ satisfies the $B$%
-condition. Note that, all uniformly convex Banach spaces with base
satisfies the $B$-condition (see $\left[ \text{10, \S\ 6},\text{ p. 75}%
\right] $, Theorem1).

\textbf{Definition 2. }A bounded linear operator\textbf{\ }$A$ is said to be
a quasi nucliar operator (QNO) of order $p$ if there is a $p\in \left(
1,\infty \right) $ such that 
\begin{equation*}
\left\Vert A\right\Vert _{p}^{p}=\left\Vert A\right\Vert _{\sigma _{p}\left(
E\right) }^{p}=\sum\limits_{i=1}^{\infty }\sum\limits_{j=1}^{\infty
}\left\vert <Ae_{i},f_{j}>\right\vert ^{p}<\infty .
\end{equation*}

The collection of such operators will be denoted by $\sigma _{p}\left(
E\right) .$

Let $s_{j}\left( A\right) $ denote the approximation numbers of the operator 
$A$ (see e.g. $\left[ \text{24, \S\ 1.16.1}\right] $). Let 
\begin{equation*}
C_{p}\left( E\right) =\left\{ A:A\in \sigma _{\infty }\left( E\right) ,\text{
}\sum\limits_{j=1}^{\infty }s_{j}^{p}\left( A\right) <\infty ,\text{ }1\leq
p<\infty \right\} .
\end{equation*}

\textbf{Remark 1. }Let $H$ be a Hilbert space and $A$ be a compact operator
in $H.$ Then $s_{j}\left( A\right) =\lambda _{j}\left( A^{\ast }A\right) ^{%
\frac{1}{2}}$, where $\lambda _{1},\lambda _{2},...$ are eigenvalues of non
negative self adjoint operator $T=\left( A^{\ast }A\right) ^{\frac{1}{2}},$
arranged in decreasing order and repeated according to multiplicity. $%
\left\{ s_{j}\left( A\right) \right\} $ are called the characteristic
numbers of the operator $A$. By Corollary 7$\ $in $\left[ \text{8, Corollary
XI. 9.1}\right] ,$ if $A\in C_{p}\left( H\right) $, $p\in \left( 0,\infty
\right) $, then the Weyl type inequality is true: 
\begin{equation}
\sum\limits_{j=1}^{\infty }\left\vert \lambda _{j}\left( A\right)
\right\vert ^{p}\leq \sum\limits_{j=1}^{\infty }s_{j}^{p}\left( A\right) .
\end{equation}

By choosing $E=H$, $A\in C_{p}\left( H\right) $ and by putting $%
f_{j}=e_{j},j=1,2,...$ in Definition 2, where $e_{j}$ are orthonormal
eigenvectors of the operator $A,$ by $\left( 1\right) $ we obtain 
\begin{equation*}
\left\Vert A\right\Vert _{\sigma _{p}\left( H\right)
}^{p}=\sum\limits_{i=1}^{\infty }\sum\limits_{j=1}^{\infty }\left\vert
\left( Ae_{i},e_{j}\right) \right\vert ^{p}=\sum\limits_{j=1}^{\infty
}\left\vert \lambda _{j}\left( A\right) \right\vert ^{p}\leq
\sum\limits_{j=1}^{\infty }s_{j}^{p}\left( A\right) =\left\Vert A\right\Vert
_{C_{p}\left( H\right) }^{p}<\infty .
\end{equation*}

It implies that $C_{p}\left( H\right) \subset \sigma _{p}\left( H\right) $.
The embedding $C_{p}\left( E\right) \subset \sigma _{p}\left( E\right) $ can
also be shown for the Banach spaces $E$ satisfying the $B$-condition. Thus,
let $E$ be a Banach space satisfying the $B$-condition and $A\in $ $%
C_{p}\left( E\right) $ such that $\left\{ e_{j}\right\} $\ is a eigen system
of the operator $A$ corresponding to the eigen values $\left\{ \lambda
_{j}\right\} $ of the $A$. So, for the appropriate biorthonormal system $%
\left\{ e_{j},f_{j}\right\} ,$ $j=1,2,...$ in $E\times $ $E^{\ast }$ we get%
\begin{equation*}
\left\vert <Ae_{i},f_{j}>\right\vert =\left\vert <\lambda
_{i}e_{i},f_{j}>\right\vert =\left\vert \lambda _{i}\right\vert .
\end{equation*}%
Then,\ by virtue of Weyl type inequality in Banach spaces $\left[ \text{12,
p. 85}\right] $ we have 
\begin{equation*}
\left\Vert A\right\Vert _{\sigma _{p}\left( E\right)
}^{p}=\sum\limits_{i=1}^{\infty }\sum\limits_{j=1}^{\infty }\left\vert
<Ae_{i},f_{j}>\right\vert ^{p}=\sum\limits_{i=1}^{\infty }\left\vert \lambda
_{i}\left( A\right) \right\vert ^{p}\leq \sum\limits_{i=1}^{\infty
}s_{i}^{p}\left( A\right) =\left\Vert A\right\Vert _{C_{p}\left( E\right)
}^{p}<\infty .
\end{equation*}%
Since all $A\in $ $C_{p}\left( E\right) $ can be approximated by sequences
of finite dimensional operators in the Banach spaces $E$ with basis$,$ the
embedding $C_{p}\left( E\right) \subset \sigma _{p}\left( E\right) $ is
shown for all $A\in C_{p}\left( E\right) .$

Let us firstly, point out some properties of the set $\sigma _{p}\left(
E\right) .$

\textbf{Corollary 1. }Let $E$ be a Banach space\ satisfying the $B$%
-condition and $A\in \sigma _{p}\left( E\right) $ for a $p\in \left(
1,\infty \right) $. Suppose $\left\{ e_{j},f_{j}\right\} ,$ $j=1,2,...$ is a
biorthonormal basis system in $E\times $ $E^{\ast }$, then there is a
positive constant $C$ so that 
\begin{equation*}
\left\Vert A\right\Vert _{p}\geq C\left( \sum\limits_{i=1}^{\infty
}\left\Vert Ae_{i}\right\Vert _{E}^{p}\right) ^{\frac{1}{p}}.
\end{equation*}

\textbf{Proof. }Really, by virtue of B-condition we have 
\begin{equation*}
\left\Vert A\right\Vert _{p}^{p}=\sum\limits_{i=1}^{\infty
}\sum\limits_{j=1}^{\infty }\left\vert <Ae_{i},f_{j}\right\vert ^{p}\geq
C\sum\limits_{i=1}^{\infty }\left\Vert Ae_{i}\right\Vert _{E}^{p}.
\end{equation*}

It is implies the assertion.

$U$ is an unitary operator in $E$ if $U$ and $U^{-1}$ are bounded in $E$ and 
$\left\Vert Ux\right\Vert _{E}=\left\Vert x\right\Vert _{E},$ $\left\Vert
U^{\ast }g\right\Vert _{E^{\ast }}=\left\Vert g\right\Vert _{E^{\ast }}$ for
all $x\in E$ and $g\in E^{\ast }$. Moreover if $\left\{ e_{j},f_{j}\right\}
, $ $j=1,2,...$ is a biorthonormal basis system in $E$ $\times E^{\ast }$,
then $\left\{ Ue_{j},\left( U^{-1}\right) ^{\ast }f_{j}\right\} $ and $%
\left\{ U^{-1}e_{j},U^{\ast }f_{i}\right\} $ are also biorthonormal basis
systems in $E$ $\times E^{\ast }.$

\textbf{Lemma 1. }Let $E$ be a Banach space satisfying the $B$-condition.
The $\sigma _{p}\left( E\right) $ norms, for a fixed $p\in \left( 1,\infty
\right) $ with respect to the different biorthonormal basis systems used in
its definition, are equivalent. If $A\in \sigma _{p}\left( E\right) $ and $U$
is a unitary operator in $E,$ then $U^{-1}AU\in \sigma _{p}\left( E\right) $
and there are positive constants $C$, $C_{1}$ and $C_{2}$ such that:

(a)%
\begin{equation*}
\left\Vert A\right\Vert _{B\left( E\right) }\leq C\left\Vert A\right\Vert
_{\sigma _{p}\left( E\right) },\text{ }\left\Vert A\right\Vert _{\sigma
_{p}\left( E\right) }=\left\Vert A^{\ast }\right\Vert _{\sigma _{p}\left(
E^{\ast }\right) }.
\end{equation*}

(b)%
\begin{equation*}
C_{1}\left\Vert A\right\Vert _{\sigma _{p}\left( E\right) }\leq \left\Vert
U^{-1}AU\right\Vert _{\sigma _{p}\left( E\right) }\leq C_{2}\left\Vert
A\right\Vert _{\sigma _{p}\left( E\right) }.
\end{equation*}

\textbf{Proof.} Suppose $\left\{ e_{j},f_{j}\right\} $ and $\left\{ \upsilon
_{j},g_{j}\right\} $, $j=1,2,...$ are two biorthonormal basis systems in $%
E\times E^{\ast }$. Then there is a unitary operator $U$ such that $\upsilon
_{j}=Ue_{j}$ and $g_{i}=\left( U^{-1}\right) ^{\ast }f_{i}.$ I.e, there are
a system of numbers $\left\{ a_{jk}\right\} ,$ $\left\{ b_{ik}\right\} $
such that $e_{j}=\sum\limits_{k=1}^{\infty }a_{jk}\upsilon _{k}$, and $%
f_{i}=\sum\limits_{m=1}^{\infty }b_{im}g_{m},$ where 
\begin{equation*}
a_{jk}=<e_{j},g_{k}>\text{and }b_{ik}=<f_{i},\upsilon _{k}>.
\end{equation*}%
Let $\left\Vert A\right\Vert _{1,p}$ and $\left\Vert A\right\Vert _{2,p}$
denote $\sigma _{p}$ norms of the operator $A$ with respect to first and
second basis systems, respectively. Substituting the above equality in the
expression $\left\Vert A\right\Vert _{1,p}^{p}$ and by using the linearity
properties of $A$ and $f_{i}$ we have 
\begin{equation*}
\left\Vert A\right\Vert _{1,p}^{p}=\sum\limits_{j=1}^{\infty
}\sum\limits_{i=1}^{\infty }\left\vert <Ae_{j},f_{i}>\right\vert
^{p}=\sum\limits_{j=1}^{\infty }\sum\limits_{i=1}^{\infty
}\sum\limits_{k=1}^{\infty }\left\vert a_{jk}\right\vert
^{p}\sum\limits_{m=1}^{\infty }\left\vert b_{im}\right\vert ^{p}\left\vert
<A\upsilon _{k},g_{m}>\right\vert ^{p}.
\end{equation*}%
By virtue of \ $B$-condition, $\sum\limits_{j=1}^{\infty }\left\vert
a_{jk}\right\vert ^{p}\leq C,$ for all $k$ and $\sum\limits_{i=1}^{\infty
}\left\vert b_{im}\right\vert ^{p}\leq C$ for all $m.$ Then we get from the
above%
\begin{equation*}
\left\Vert A\right\Vert _{1,p}^{p}\leq C_{1}\sum\limits_{k=1}^{\infty
}\sum\limits_{m=1}^{\infty }\left\vert <A\upsilon _{k},g_{m}>\right\vert
^{p}=C_{1}\left\Vert A\right\Vert _{2,p}^{p}.
\end{equation*}

In a similar way, we get 
\begin{equation*}
\left\Vert A\right\Vert _{2,p}^{p}\leq C_{2}\left\Vert A\right\Vert
_{1,p}^{p}.
\end{equation*}

This implies that $\sigma _{p}\left( E\right) $ norms are independent of the
biorthonormal basis systems.

Let $\left\{ e_{j},f_{j}\right\} ,$ $j=1,2,...$ be a biorthonormal basis
system in $E\times E^{\ast }$. By using Definition 2 it is seen that 
\begin{equation*}
\left\Vert A\right\Vert _{p}^{p}=\sum\limits_{i=1}^{\infty
}\sum\limits_{j=1}^{\infty }\left\vert <Ae_{i},f_{j}>\right\vert
^{p}=\sum\limits_{i=1}^{\infty }\sum\limits_{j=1}^{\infty }\left\vert
<e_{i},A^{\ast }f_{j}>\right\vert ^{p}=\left\Vert A^{\ast }\right\Vert
_{p}^{p}.
\end{equation*}%
The assertion ( b) is obtained from the equivalence of $\sigma _{p}\left(
E\right) $ norms with respect to different basis systems. Really, if $U$ is
a uniter operator in $E$, then $\left\{ Ue_{i},\left( U^{-1}\right) ^{\ast
}f_{j}\right\} $ is a biorthonormal system in $E\times E^{\ast }$. So, we
have%
\begin{equation*}
\left\Vert U^{-1}AU\right\Vert _{p}^{p}=\sum\limits_{i=1}^{\infty
}\sum\limits_{j=1}^{\infty }\left\vert <U^{-1}AUe_{i},f_{j}>\right\vert ^{p}
\end{equation*}

\begin{equation*}
=\sum\limits_{i=1}^{\infty }\sum\limits_{j=1}^{\infty }\left\vert
<AUe_{i},\left( U^{-1}\right) ^{\ast }f_{j}>\right\vert ^{p}\leq
C_{2}\left\Vert A\right\Vert _{p}^{p}.
\end{equation*}%
In a similar way we get%
\begin{equation*}
\left\Vert A\right\Vert _{p}^{p}\leq C_{1}\left\Vert U^{-1}AU\right\Vert
_{p}^{p}.
\end{equation*}
These two inequalities imply the assumption (b).

Finally, if $\varepsilon >0$ let $u_{0}$ be an element of unit norm such
that 
\begin{equation*}
\left\Vert A\right\Vert _{B\left( E\right) }^{p}\leq \left\Vert
Au_{0}\right\Vert _{E}^{p}+\varepsilon .
\end{equation*}

Then, by definition of $\sigma _{p}\left( E\right) $ and by Corallary1 we
get 
\begin{equation*}
\left\Vert A\right\Vert _{B\left( E\right) }\leq C\left\Vert A\right\Vert
_{\sigma _{p}\left( E\right) }.
\end{equation*}

\textbf{Remark 2. }The basis equivalence of $\sigma _{p}\left( E\right) $
norms, for a fixed $p\in \left( 1,\infty \right) ,$ mean that, there are the
positive constants $C_{1},$ $C_{2}$ such that $\left\Vert A\right\Vert
_{i,p},$ $i=1,2$ norms with respect to different two biorthonormal basis
systems satisfy the relation 
\begin{equation*}
C_{1}\left\Vert A\right\Vert _{1,p}^{p}\leq \left\Vert A\right\Vert
_{2,p}^{p}\leq C_{2}\left\Vert A\right\Vert _{1,p}^{p}.
\end{equation*}%
The independence of $\sigma _{p}\left( E\right) $ norms of basis systems are
valid when $E$ is a Hilbert space.

In a similar way as in $\left[ \text{8, Theorem XI. 6.4-7, }\right] $ we have

\textbf{Theorem A}$_{1}$\textbf{. }Let $E$ be a Banach space satisfying the $%
B$-condition. Then, the set $\sigma _{p}\left( E\right) ,$ $p\in \left(
1,\infty \right) $ is a Banach space under $\sigma _{p}\left( E\right) $
norm.

\textbf{Theorem A}$_{2}$\textbf{. }Let $E$ be a Banach space satisfying the $%
B$-condition. Then, every $A\in \sigma _{p}\left( E\right) $, $p\in \left(
1,\infty \right) $ is a compact operator in $E$ and is a limit in $\sigma
_{p}\left( E\right) $ norm of a sequence of operators with finite
dimensional range.

\textbf{Theorem A}$_{3}$\textbf{. }Let $E$ be a Banach space satisfying the $%
B$-condition. If $A\in \sigma _{p}\left( E\right) $ for a $p\in \left(
1,\infty \right) $ and $F$ is a single-valued analytic function on its
spectrum which vanishes at zero, then $F\left( A\right) $ $\in \sigma
_{p}\left( E\right) $ and the map $A\rightarrow F\left( A\right) $ is
continuous in $\sigma _{p}\left( E\right) $. Furthermore, if $\left\{
F_{n}\right\} $ is a sequence of such functions having as common domain a
neighborhood $N$ of the spectrum of $A$ and if $F_{n}\left( \lambda \right)
\rightarrow F\left( \lambda \right) $ uniformly for $\lambda $ in $N$, then $%
F_{n}\left( A\right) \rightarrow F\left( A\right) $ in $\sigma _{p}\left(
E\right) .$

\textbf{Lemma 2. }Let $E$ be a Banach space satisfying the $B$-condition and 
$A\in \sigma _{p}\left( E\right) ,$ $B\in \sigma _{q}\left( E\right) $ for a 
$p,q\in \left( 1,\infty \right) ,$ where $\frac{1}{p}+\frac{1}{q}=1.$
Suppose $\left\{ e_{j},f_{j}\right\} ,$ $j=1,2,...$ is a biorthonormal basis
system in $E\times E^{\ast }$, then the series $\sum\limits_{i=1}^{\infty
}<Ae_{i},B^{\ast }f_{i}>$ converges absolutely to a limit which is
independent of the basis. Moreover, 
\begin{equation*}
\sum\limits_{i=1}^{\infty }<Ae_{i},B^{\ast }f_{i}>=\sum\limits_{i=1}^{\infty
}<Be_{i},A^{\ast }f_{i}>.
\end{equation*}

\textbf{Proof. }By Holder inequality we have 
\begin{equation*}
\sum\limits_{i=1}^{\infty }\sum\limits_{j=1}^{\infty }\left\vert
<Ae_{i},f_{j}><e_{i},B^{\ast }f_{j}>\right\vert \leq \left\{
\sum\limits_{i=1}^{\infty }\sum\limits_{j=1}^{\infty }\left\vert
<Ae_{i},f_{j}>\right\vert ^{p}\right\} ^{\frac{1}{p}}
\end{equation*}%
\begin{equation*}
\left\{ \sum\limits_{i=1}^{\infty }\sum\limits_{j=1}^{\infty }\left\vert
<e_{i},B^{\ast }f_{j}>\right\vert ^{q}\right\} ^{\frac{1}{q}}=\left\Vert
A\right\Vert _{\sigma _{p}\left( E\right) }\left\Vert B^{\ast }\right\Vert
_{\sigma _{q}\left( E\right) }.
\end{equation*}

Thus the double series $\sum\limits_{i=1}^{\infty }\sum\limits_{j=1}^{\infty
}<Ae_{i},f_{j}><e_{i},B^{\ast }f_{j}>$ converges absolutely, and hence the
corresponding iterated series exists and are equal. Moreover, by $B$%
-condition, there is another biorthonormal basis system $\left\{
e_{j}^{^{\prime }},f_{j}^{^{\prime }}\right\} ,$ $j=1,$ $2,...$ in $E\times
E^{\ast }$ such that 
\begin{equation}
\sum\limits_{i=1}^{\infty }<Ae_{i},B^{\ast }f_{i}>=\sum\limits_{i=1}^{\infty
}\sum\limits_{j=1}^{\infty }<Ae_{i},f_{j}^{^{\prime }}><e_{i}^{\prime
},B^{\ast }f_{j}>
\end{equation}%
\begin{equation*}
=\sum\limits_{j=1}^{\infty }\sum\limits_{i=1}^{\infty }<e_{i},A^{\ast
}f_{j}^{^{\prime }}><Be_{i}^{^{\prime }},f_{j}>=\sum\limits_{j=1}^{\infty
}<Be_{i}^{\prime },A^{\ast }f_{j}^{^{\prime }}>.
\end{equation*}

From $\left( 2\right) $ we obtain 
\begin{equation*}
\sum\limits_{i=1}^{\infty }<Ae_{i},B^{\ast }f_{i}>=\sum\limits_{i=1}^{\infty
}<Be_{i},A^{\ast }f_{i}>.
\end{equation*}

Hence, this expression is symmetric in $A$ and $B.$ By using $\left(
2\right) $ we obtain the independence of the limit from the basis systems.

\textbf{Definition 3. }Let $E$ be a Banach space\ satisfying the $B$%
-condition and 
\begin{equation*}
A\in \sigma _{p}\left( E\right) ,\text{ }B\in \sigma _{q}\left( E\right) 
\text{ for a }p,\text{ }q\in \left( 1,\infty \right) ,\frac{1}{p}+\frac{1}{q}%
=1.
\end{equation*}%
Suppose $\left\{ e_{i},f_{i}\right\} ,$ $i=1,2,...$ is a biorthonormal basis
system in $E\times E^{\ast }$, then the trace of $\left( A,B\right) $ is
defined to be as: 
\begin{equation*}
Tr\left( A,B\right) =\sum\limits_{i=1}^{\infty }<Ae_{i},B^{\ast }f_{i}>.
\end{equation*}

\textbf{Corollary 2.} Let $A\in \sigma _{p}\left( E\right) ,$ $B\in \sigma
_{q}\left( E\right) $ for a $p,$ $q\in \left( 1,\infty \right) $, $\frac{1}{p%
}+\frac{1}{q}=1,$ then the trace is a symmetric bilinear function and 
\begin{equation}
Tr\left( A,B\right) \leq \left\Vert A\right\Vert _{\sigma _{p}\left(
E\right) }\left\Vert B^{\ast }\right\Vert _{\sigma _{q}\left( E\right) }.
\end{equation}

\textbf{Proof. }The symmetry of the trace function were proved during the
proof of Lemma 2. Moreover, by $\left( 2\right) $ we get 
\begin{equation*}
\sum\limits_{i=1}^{\infty }<Ae_{i},B^{\ast }f_{i}>=\sum\limits_{i=1}^{\infty
}\sum\limits_{j=1}^{\infty }<Ae_{i},f_{j}^{^{\prime }}><e_{j}^{\prime
},B^{\ast }f_{j}>.
\end{equation*}%
So by Holder inequality and B-condition we have 
\begin{equation*}
\sum\limits_{i=1}^{\infty }<Ae_{i},B^{\ast }f_{i}>\leq \left\{
\sum\limits_{i=1}^{\infty }\sum\limits_{j=1}^{\infty }\left\vert
<Ae_{i},f_{j}^{^{\prime }}>\right\vert ^{p}\right\} ^{\frac{1}{p}}
\end{equation*}%
\begin{equation*}
\left\{ \sum\limits_{i=1}^{\infty }\sum\limits_{j=1}^{\infty }\left\vert
<e_{i},B^{\ast }f_{j}>\right\vert ^{q}\right\} ^{\frac{1}{q}}=\left\Vert
A\right\Vert _{\sigma _{p}\left( E\right) }\left\Vert B^{\ast }\right\Vert
_{\sigma _{q}\left( E\right) }.
\end{equation*}

This relation implies the assertion$.$

In a similar manner as $\left[ \text{8, Lemma XI. 6.20 }\right] $ we have

\textbf{Lemma A}$_{1}.$ Let $E$ be a Banach space satisfying the
B-condition. Suppose $A\in \sigma _{p}\left( E\right) $ for a $p\in \left(
1,\infty \right) $ having a finite dimensional range. Let $N\left( A\right) $
be the null space of $A,$ and let $P$ be the orthogonal projection onto a
finite dimensional subspace of $E$ containing $\left[ N\left( A\right) %
\right] ^{\perp }.$ Then:

(a) the spectra of the operators $A$ and $PA$ coincide;

(b) For a single valued analytic function $F$ on spectrum of $A$ with $%
F\left( 0\right) =0$, the following hold 
\begin{equation*}
F\left( PA\right) =PF\left( A\right) ,F\left( A\right) =F\left( A\right) P%
\text{;}
\end{equation*}

(c) $Tr\left( F\left( A\right) ,A\right) =Tr\left( F\left( PA\right)
,PA\right) $ and $Tr\left( F\left( PA\right) ,PA\right) $ coincide with the
trace of the restriction of the operator $PAF\left( A\right) $ to the finite
dimensional space $PE.$

\textbf{Proof. }The (a) and (b) parts are proving by using the spectral
properties of compact operators and operator calculus as in $\left[ \text{8,
Lemma XI. 6.20 }\right] .$ Let $\left\{ e_{j},f_{j}\right\} ,$ $j=1,2,...$
is a biorthonormal basis system in $E\times E^{\ast }$. Since $PE$ is finite
dimensional we may suppose that there is a number $d$ such that finite set $%
\left\{ e_{j}\right\} ,j=1,2,..,d$ is a basis for $PE,$ and the sub set $%
\left\{ e_{j}\right\} ,$ $j=d+1,d+2,...$ is a basis for $\left( I-P\right)
E. $ Then, since $A=AP,$ we have $A^{\ast }=PA^{\ast }$ and 
\begin{equation*}
Tr\left( F\left( PA\right) ,PA\right) =Tr\left( PF\left( A\right) ,PA\right)
=\sum\limits_{j=1}^{\infty }<PF\left( A\right) e_{j},\left( PA\right) ^{\ast
}f_{j}>
\end{equation*}%
\begin{equation*}
=\sum\limits_{j=1}^{d}<F\left( A\right) e_{j},A^{\ast
}f_{j}>=\sum\limits_{j=1}^{d}<PAF\left( A\right) e_{j},f_{j}>=Tr\left[
PAF\left( A\right) \mid PE\right] .
\end{equation*}

Since $F\left( A\right) E=F\left( A\right) ,$ we have $F\left( A\right) $ $%
\left( I-P\right) =0,$ $F\left( A\right) e_{j}=0$ for $j=d+1,d+2,...,$ and
it follows from the above that 
\begin{equation*}
Tr\left( F\left( A\right) ,A\right) =\sum\limits_{j=1}^{\infty }<F\left(
A\right) e_{j},A^{\ast }f_{j}>=\sum\limits_{j=1}^{d}<F\left( A\right)
e_{j},A^{\ast }f_{j}>
\end{equation*}

\begin{equation*}
=Tr\left( F\left( PA\right) ,PA\right)
\end{equation*}%
which implies the (c) part.

In a similar way as $\left[ \text{8, Lemmas XI. 6.21- 6.23 }\right] $ we
obtain, respectively.

\textbf{Lemma A}$_{2}.$ Let $\lambda $ and $z$ be complex numbers with $%
\lambda z\neq 1$ and let 
\begin{equation*}
F\left( \lambda ,z\right) =z^{-1}\left[ \log \left( 1-\lambda z\right)
+\lambda z\right] .
\end{equation*}

Let $E$ be a Banach space\ satisfying the $B$-condition and $A\in \sigma
_{p}\left( E\right) $ for a $p\in \left( 1,\infty \right) $ whose spectrum
does not include the number $\lambda ^{-1}$. Suppose $\left\{
e_{j},f_{j}\right\} ,$ $j=1,2,...$ is a biorthonormal basis system in $%
E\times E^{\ast }$. Then for any finite subsets $\left\{ e_{j}\right\}
,\left\{ f_{j}\right\} ,$ $j=1,2,..,d$ \ the following inequality holds: 
\begin{equation*}
\exp \left[ Tr\left( F\left( \lambda ,A\right) ,A\right) \right] \leq \exp
\left\{ \frac{1}{p}\sum\limits_{j=1}^{d}\left\vert \lambda Ae_{j}\right\vert
^{p}\right\}
\end{equation*}%
\begin{equation*}
\exp \left\{ \frac{1}{p}\sum\limits_{j=1}^{d}\func{Re}<\lambda
Ae_{j},f_{j}>\right\} \prod\limits_{j=1}^{d}\left[ 1-2\func{Re}\left(
<\lambda Ae_{j},f_{j}>\right) +\left\Vert \lambda Ae_{j}\right\Vert ^{p}%
\right] ^{\frac{1}{p}}.
\end{equation*}

\textbf{Lemma A}$_{3}.$ For any positive $\varepsilon $ we have 
\begin{equation*}
\lim\limits_{\left\vert \lambda \right\vert \rightarrow \infty
}e^{-\varepsilon \left\vert \lambda \right\vert ^{2}}\exp \left[ Tr\left(
F\left( \lambda ,A\right) ,A\right) \right] =0.
\end{equation*}

In a similar way as $\left[ \text{8, Theorem XI. 6.24}\right] $ we have

\textbf{Theorem A}$_{4}$\textbf{.} Let $E$ be a Banach space satisfying the
B-condition. Assume $N\in \sigma _{p}\left( E\right) $ for a $p\in \left(
1,\infty \right) $ is a quasi-nilpotent operator. Then $Tr\left( N,N\right)
=0.$

We are now in a position to obtain results in infinite dimensional Banach
spaces by using of key finite dimensional results. By this aim by following $%
\left[ \text{8, Theorem XI. 6.24}\right] $\ we obtain

\textbf{Theorem 1. }Let $E$ be a Banach space\ satisfying the $B$-condition.
Suppose$\ A\in \sigma _{p}\left( E\right) $ for a $p\in \left( 1,\infty
\right) $ and $\lambda _{1},\lambda _{2},...$ are its eigenvalues repeated
according to multiplicities. If $F$ and $g$ are functions analytic in a
neighborhood of the spectrum of $A$ with $F\left( 0\right) =0,$ $g\left(
0\right) =0,$ then $F\left( A\right) $, $g\left( A\right) \in \sigma
_{p}\left( E\right) ,$ and 
\begin{equation*}
Tr\left( F\left( A\right) ,g\left( A\right) \right)
=\sum\limits_{i=1}^{\infty }F\left( \lambda _{i}\right) g\left( \lambda
_{i}\right) ,
\end{equation*}%
where the series on the right hand side is absolutely convergent.

\textbf{Proof.} At first, by reasoning as the beginning of the proof $\left[ 
\text{8, Theorem XI. 6.25}\right] ,$ we get 
\begin{equation*}
\sum\limits_{i=1}^{\infty }\left\vert F\left( \lambda _{i}\right)
\right\vert ^{p}<\infty ,\text{ }\sum\limits_{i=1}^{\infty }\left\vert
F\left( \lambda _{i}\right) g\left( \lambda _{i}\right) \right\vert <\infty .
\end{equation*}%
Let $P_{i}=P\left( \lambda _{i};A\right) $ denote the projection operators
defined in $\left[ \text{8, }\vee 11\text{.3}\right] $ i.e. 
\begin{equation*}
P_{i}E=E_{i},\text{ dim }E_{i}<\infty \text{, }i=1,2,....
\end{equation*}%
Let $G_{1}$ be the closure of the subspace $\sum\limits_{i=1}^{\infty
}P_{i}E $ and $G_{2}$ be the orthocomplement of the $G_{1}$, i.e. 
\begin{equation*}
G_{2}=\left\{ f\in E^{\ast }:<u,f>=0,\text{ }u\in G_{1}\right\} .
\end{equation*}%
Suppose $\left\{ e_{j},f_{j}\right\} ,$ $j=1,2,...$ is a biorthonormal basis
system in $E\times E^{\ast }$. Assume $\left\{ e_{j}\right\} $ so that the
sub system$\left\{ e_{j}\right\} $, $j=1,2,...n_{1}$ is a basis for $E_{1},$ 
$\left\{ e_{j}\right\} $, $j=1,2,...n_{2}$ is a basis for $E_{2}$, etc. Let $%
\left\{ \Psi _{k}\right\} $ be a sub system of $\left\{ f_{j}\right\}
\subset E^{\ast }$ which is a basis for $G_{2}.$ Then by Definition 3 and
Theorem A$_{3}$ we get 
\begin{equation*}
Tr\left( F\left( A\right) ,g\left( A\right) \right)
=\sum\limits_{j=1}^{\infty }<F\left( A\right) e_{j},\left( g\left( A\right)
\right) ^{\ast }f_{j}>+\sum\limits_{k=1}^{\infty }<F\left( A\right)
e_{k},\left( g\left( A\right) \right) ^{\ast }\Psi _{k}>.
\end{equation*}%
By Theorems A$_{1}$-A$_{2}$ and Lemma A$_{1}$ we have 
\begin{equation*}
\sum\limits_{j=1}^{\infty }<F\left( A\right) e_{j},\left( g\left( A\right)
\right) ^{\ast }f_{j}>=\lim\limits_{j\rightarrow \infty
}\sum\limits_{j=1}^{n_{j}}<F\left( A\right) e_{j},\left( g\left( A\right)
\right) ^{\ast }f_{j}>
\end{equation*}%
\begin{equation*}
=\lim\limits_{j\rightarrow \infty }Tr\left( gF\left( A\right) ,AP_{j}\right)
=\sum\limits_{i=1}^{\infty }F\left( \lambda _{i}\right) g\left( \lambda
_{i}\right) .
\end{equation*}

Now it is sufficient to show the equality 
\begin{equation}
\sum\limits_{k=1}^{\infty }<F\left( A\right) e_{k},\left( g\left( A\right)
\right) ^{\ast }\Psi _{k}>=0.
\end{equation}

By Lemma 2 we have 
\begin{equation*}
<F\left( A\right) e_{k},\left( g\left( A\right) \right) ^{\ast }\Psi
_{k}>=<g\left( A\right) e_{k},\left( F\left( A\right) \right) ^{\ast }\Psi
_{k}>.
\end{equation*}

So, the validity of $\left( 3\right) $ is a consequence of the validity of
the following equations 
\begin{equation*}
\sum\limits_{k=1}^{\infty }<F\left( A\right) e_{k},\left( F\left( A\right)
\right) ^{\ast }\Psi _{k}>=0,\text{ }\sum\limits_{k=1}^{\infty }<g\left(
A\right) e_{k},\left( g\left( A\right) \right) ^{\ast }\Psi _{k}>=0,
\end{equation*}%
\begin{equation}
\sum\limits_{k=1}^{\infty }<\left( F+g\right) Ae_{k},\left( F+g\right)
\left( A\right) ^{\ast }\Psi _{k}>=0.
\end{equation}

All these equations being of the same forme. So it is sufficient to show one
of them. Let us prove the first of them.

By $\left[ \text{8, Theorem}\vee 11\text{.3.20}\right] ,$ $G_{1}$ is mapped
into itself by $F\left( A\right) .$ Thus $G_{2}$ is mapped into itself by $%
F\left( A\right) ^{\ast }.$ Let 
\begin{equation*}
F\left( A\right) ^{\ast }\left\vert G_{2}\right. =S.
\end{equation*}

Then by Theorem A$_{3}$, Lemma1, and Definition 2 we get $S\in \sigma
_{p}\left( E\right) $ and 
\begin{equation*}
<PF\left( A\right) u,\upsilon >=<F\left( A\right) u,\upsilon >=<u,F\left(
A\right) ^{\ast }\upsilon >,\text{ }u,\upsilon \in G_{2},
\end{equation*}%
where $P$ denoted the projection of $E$ on $G_{2}.$ Thus $PF\left( A\right)
\mid G_{2}=S^{\ast }.$ Hence $\left( 4\right) $ is equivalent to the
assertion 
\begin{equation}
Tr\left( S,S\right) =0.
\end{equation}

It follows from Theorem A$_{4}$ that to prove $\left( 5\right) $, it
sufficient to show that $S$ is quasi-nilpotent. If this is not so, then by $%
\left[ \text{8, Theorem}\vee 11\text{.4.5}\right] ,$ there exists a non-zero
complex number $\mu $ and a non-zero element $u\in G_{2}$ such that $Su=\mu
u.$ Thus, by $\left[ \text{8, Theorem}\vee 11\text{.4.5}\right] $ again, $%
P\left( \mu ,F\left( A\right) ^{\ast }\right) G_{2}\neq \left\{ 0\right\} .$
By definitions $\left[ \text{8, }\vee 11.3.9,3.17,\right] $ and by $\left[ 
\text{8, Lemma }1\vee .2.10\right] ,$ it is seen that 
\begin{equation*}
P\left( \mu ,F\left( A\right) ^{\ast }\right) =\left( P\left( \bar{\mu}%
,F\left( A\right) \right) \right) ^{\ast }.
\end{equation*}

Hence, according to $\left[ \text{8, Theorem}\vee 11\text{.3.20}\right] $,
there is a non-zero complex number $\nu $ such that $P\left( \nu ,A^{\ast
}\right) G_{2}\neq \left\{ 0\right\} .$ However, since $<G_{2},P\left( \nu ,%
\text{ }A^{\ast }\right) E>=0$ for $\nu \neq 0,$ by definition we have a
contradiction which proves the present theorem.

In a similar way as $\left[ \text{8, Theorem XI. 6.26}\right] $ we have

\textbf{Theorem A}$_{5}$\textbf{. }Assume $E$ is a Banach space satisfying
the $B$-condition. Let $A\in \sigma _{p}\left( E\right) $ for a $p\in \left(
1,\infty \right) $ and let $\lambda _{1},\lambda _{2},...,\lambda _{n}...$
be its eigenvalues repeated according to multiplicities. Then the infinite
product $\varphi _{\lambda }\left( A\right) =\prod\limits_{i=1}^{\infty
}\left( 1-\frac{\lambda _{i}}{\lambda }\right) e^{\frac{\lambda _{i}}{%
\lambda }}$ converges and defines a function analytic for $\lambda \neq 0.$
For each fixed $\lambda \neq 0$ and $\varphi _{\lambda }\left( A\right) $ is
a continuous complex valued function on the Banach space of $\sigma
_{p}\left( E\right) .$

Now we can state the following Carleman theorem in Banach spaces.

\textbf{Theorem 2. }Let $E$ be a Banach space\ satisfying the $B$-condition.
Let $A\in \sigma _{p}\left( E\right) $ for a $p\in \left( 1,\infty \right) .$
If $\lambda $ is in the resolvent set of the operator $A$, then 
\begin{equation*}
\left\Vert \varphi _{\lambda }\left( A\right) \left( \lambda -A\right)
^{-1}\right\Vert _{B\left( E\right) }\leq \left\vert \lambda \right\vert
\exp \left\{ \frac{1}{2}\left( 1+\frac{\left\Vert A\right\Vert _{p}^{p}}{%
\left\vert \lambda \right\vert ^{2}}\right) \right\} .
\end{equation*}

\textbf{Proof. }It follows from Theorem A$_{5}$ and $\left[ \text{8, Lemma}%
\vee 11.\text{6.1}\right] $, that it is sufficient to consider the case in
which $A$ has a finite dimensional range $R\left( A\right) .$ Let $N\left(
A\right) =\left\{ u\in E:\text{ }Au=0\right\} .$ Then $E/N\left( A\right) $
is mapped by $A$ in a one-to-one fashion into $R\left( A\right) .$ Thus $%
E/N\left( A\right) $ is the finite dimensional space. Let $V$ \ be a one
dimensional subspace of $N\left( A\right) $, $V_{1}=E/N\left( A\right)
+R\left( A\right) +V$ and $V_{2}=E/V_{1}.$ Then $AV_{2}=0,$ and $%
AV_{1}\subset V_{1}.$ Put $A_{1}=A\left\vert V_{1}\right. $. Then it is easy
to see that 
\begin{equation*}
\left\Vert A_{1}\right\Vert _{\sigma _{p}\left( E\right) }=\left\Vert
A\right\Vert _{\sigma _{p}\left( E\right) },\text{ }\sigma \left(
A_{1}\right) =\sigma \left( A\right) ,\text{ }\varphi _{\lambda }\left(
A_{1}\right) =\varphi _{\lambda }\left( A\right) .
\end{equation*}

Moreover, if $u_{i}\in V_{i},$ $i=1,2$, then 
\begin{equation*}
\left( \lambda -A\right) ^{-1}\left( u_{1}+u_{2}\right) =\left( \lambda
-A_{1}\right) ^{-1}u_{1}+\lambda ^{-1}u_{2}.
\end{equation*}

Thus 
\begin{equation*}
\left\| \left( \lambda -A\right) ^{-1}\right\| =\max \left\{ \left| \lambda
^{-1}\right| ,\left\| \left( \lambda -A_{1}\right) ^{-1}\right\| \right\} .
\end{equation*}

On the other hand we have 
\begin{equation*}
\left\| \left( \lambda -A_{1}\right) ^{-1}\right\| \geq \left| \lambda
^{-1}\right| .
\end{equation*}

Really, if we suppose$\left\Vert \left( \lambda -A_{1}\right)
^{-1}\right\Vert <\left\vert \lambda ^{-1}\right\vert $, then $\left[ \text{%
8, Lemma}\vee \text{11.6.1}\right] $ imply that $A_{1}$ had an inverse which
is impossible since the eigenvectors in $V$ belong to its domain $V_{1}$.
Thus 
\begin{equation*}
\left\Vert \left( \lambda -A\right) ^{-1}\right\Vert =\left\Vert \left(
\lambda -A_{1}\right) ^{-1}\right\Vert .
\end{equation*}

Consequently, the present theorem follows immediately from $\left[ \text{8,
Theorem XI. 15}\right] .$

Theorem 2 implies the following

\textbf{Corollary 3. }Let $E$ be a Banach space satisfying the $B$%
-condition. Let $N$ be a quasi-nilpotent operator in $\sigma _{p}\left(
E\right) $ for a $p\in \left( 1,\infty \right) .$ Then for every $\lambda
\neq 0$ we have 
\begin{equation*}
\left\Vert \left( \lambda -N\right) ^{-1}\right\Vert \leq \left\vert \lambda
\right\vert \exp \left\{ M\left( 1+\left\Vert \frac{N}{\lambda }\right\Vert
_{\sigma _{p}}^{p}\right) \right\} \text{, }M>0.
\end{equation*}

Now we are a position to prove the main theorem.

\textbf{Theorem 3. }Assume:

(1) $E$ is a Banach space\ satisfying the $B$-condition and $A$ is an
operator in $\sigma _{p}\left( E\right) $ for a $\ p\in \left( 1,\infty
\right) ;$

(2) $\gamma _{1},\gamma _{2},...,\gamma _{s}$ is non overlapping,
differentiable arcs in the complex plane starting at the origin. Suppose
that each of the $s$ regions into which the plans is divided by these arcs
is contained in an angular sector of opening less then $\frac{\pi }{p};$

(3) $m>0$ is an integer so that the resolvent of $A$ satisfies the
inequality $\left\Vert R\left( \lambda ,A\right) \right\Vert =O\left(
\left\vert \lambda \right\vert ^{-m}\right) $ \ as $\lambda \rightarrow 0$
along any of the arcs $\gamma _{i\text{ }}.$

Then the subspace $spA$ contains the subspace $A^{m}E.$

\textbf{Proof. }By the Hahn-Banach theorem it suffices to prove that every
element $f\in E^{\ast }$ satisfying the condition $<u,f>=0$ for $u\in spA$
also has $<A^{m}u,f>=0$ for all $u\in E.$ Let $f$ \ be such element. By
theorem $\left[ \text{8, Theorem}\vee 11.\text{4.5}\right] ,$ the function $%
f\left( \lambda \right) =\lambda ^{m}R\left( \lambda ,A^{\ast }\right) f$ \
is analytic everywhere in the plane except at $\lambda =0$ and at an
isolated set of points $\lambda _{k}\rightarrow \infty $, and at the points $%
\lambda _{k}$ the function $f\left( \lambda \right) $ may have a pole. For $%
\lambda \neq \lambda _{k}$ and $\lambda $ in the neighborhood of $\lambda
_{k}$ we have 
\begin{equation*}
f\left( \lambda \right) =\lambda ^{m}P\left( \lambda _{k},A^{\ast }\right)
R\left( \lambda ,A^{\ast }\right) f+\lambda ^{m}R\left( \lambda ,A^{\ast
}\right) \left( I-P\left( \lambda _{k},A^{\ast }\right) \right) f=
\end{equation*}%
\begin{equation*}
\lambda ^{m}P\left( \bar{\lambda}_{k},A\right) ^{\ast }R\left( \bar{\lambda}%
,A\right) ^{\ast }f+f_{1}\left( \lambda \right) .
\end{equation*}

By virtue of $\left[ \text{8, Theorem }\vee 11.\text{3.20}\right] $ and $%
\left[ \text{8, Lemma}\vee 11.\text{3.2}\right] $ the function $f_{1}\left(
\lambda \right) $ is analytic at $\lambda =\lambda _{k}.$ It will now be
shown that the function $f_{2}\left( \lambda \right) =\lambda ^{m}P\left( 
\bar{\lambda}_{k},A\right) ^{\ast }R\left( \bar{\lambda},A\right) ^{\ast }f$
\ vanishes which will prove that $f\left( \lambda \right) $ is analytic at
all the points $\lambda =\lambda _{k},$ so that $f\left( \lambda \right) $
can only fail to be analytic at the point $\lambda =0.$ Really, note that 
\begin{equation}
<u,f_{2}\left( \lambda \right) >=<u,\lambda ^{m}P\left( \bar{\lambda}%
_{k},A\right) ^{\ast }R\left( \bar{\lambda},A\right) ^{\ast }f>
\end{equation}

\begin{equation*}
=\lambda ^{m}<P\left( \bar{\lambda}_{k},A\right) R\left( \bar{\lambda}%
,A\right) u,f>.
\end{equation*}%
It follows from $\left[ \text{8, Theorem }\vee 11.\text{4.5}\right] $ that 
\begin{equation*}
P\left( \bar{\lambda}_{k},A\right) R\left( \bar{\lambda},A\right) u\in spA.
\end{equation*}

Since $f\in \left( spA\right) ^{\perp }$, $\left( 6\right) $ implies that $%
<u,f_{2}\left( \lambda \right) >=0$ for every $u\in E$ and thus $f_{2}\left(
\lambda \right) =0.$ Therefore $\lambda ^{m}R\left( \lambda ,A^{\ast
}\right) f$ \ is analytic everywhere in the plane except at the origin. If
the function $f$ is analytic at the origin then by reasoning as in $\left[ 
\text{8, Theorem XI. 6.29}\right] $ and by Liouville's theorem we obtain the
assertion$.$ So the proof rests upon the assertion that the function $%
f\left( \lambda \right) $ is analytic at $\lambda =0.$ By using the
Corollary 3, in a similar way as $\left[ \text{8, Theorem XI. 6.29}\right] ,$
we get that 
\begin{equation*}
\left\Vert R\left( \lambda ,A\right) \right\Vert =O\left( \exp \left\{
M\left( 1+\left\Vert \frac{N}{\lambda }\right\Vert _{\sigma _{p}}^{p}\right)
\right\} \right) ,\text{ }M>0
\end{equation*}

as $\lambda \rightarrow 0$. Then by virtue of Phragmen-Lindel\"{o}f theorem
we obtain that the function $f$ \ is analytic at the origin.

By using Theorem 3, in a similar way as $\left[ \text{8, Corollary XI. 6.30}%
\right] $ we have

\textbf{Corollary 4.} Suppose (1) and (2) condition of Theorem 3 hold and
resolvent of $A$ satisfies the inequality $\left\Vert R\left( \lambda \text{,%
}A\right) \right\Vert =O\left( \left\vert \lambda \right\vert ^{-1}\right) $
\ as $\lambda \rightarrow \infty $ along any of the arcs $\gamma _{i\text{ }%
}.$ Then the subspace $spA$ contains the subspace $E.$

\textbf{Proof: }It is sufficient to show that joint span of the range $%
R\left( A\right) $ and the null space $N\left( A\right) $ is the entire
space $E.$ Let $\left\{ \lambda _{n}\right\} $ be a sequence of complex
numbers converging to zero along one of the arcs $\gamma _{i}$ and let $u$
be an arbitrary element from $E.$ By assumptions, the sequence $\left\{
\lambda _{n}R\left( \lambda _{n}\text{, }A\right) \right\} $ is bounded.
Since $E$ is reflexive, then this sequence is weakly convergent to an
element $\upsilon .$ The proof will be completed by showing that $A\upsilon
=0$ and $u-\upsilon \in \bar{N}\left( A\right) .$ Then, by reasoning as in
the proof of $\left[ \text{8, Corollary XI. 6.30}\right] $ we obtain the
assertion.

By using Theorem 3, in a similar way as $\left[ \text{8, Corollary XI. 6.31}%
\right] $ we have

\textbf{Corollary 5. }Suppose:

(1) $E$ is a Banach space\ satisfying the $B$-condition;

(2) $A$ is a densely defined unbounded operator in $E,$ with the property
that for some $\lambda $ in the resolvent, the operator $R\left( \lambda 
\text{, }A\right) $ is of class $\sigma _{p}\left( E\right) $ for a $p\in
\left( 1,\infty \right) ;$

(2) $\gamma _{1},\gamma _{2},...,\gamma _{s}$ is non overlapping,
differentiable arcs in the complex plane having a limiting direction at
infinity, and such that no adjacent pair of arcs form an angle as great as $%
\frac{\pi }{p}$ at infinity;

(3) the resolvent of $A$ satisfies the inequality $\left\Vert R\left(
\lambda ,A\right) \right\Vert =O\left( \left\vert \lambda \right\vert
^{-1}\right) $ \ as $\lambda \rightarrow \infty $ along any of the arcs $%
\gamma _{i\text{ }}.$

Then the subspace $spA$ contains the entire space $E.$

\begin{center}
\textbf{Spectral properties of abstract elliptic operators}
\end{center}

Consider the nonlocal BVP for differential operator equation%
\begin{equation}
\left( L+\lambda \right) u=a\left( x\right) u^{\left( 2\right) }\left(
x\right) +B\left( x\right) u^{^{\left( 1\right) }}\left( x\right)
+A_{\lambda }\left( x\right) u\left( x\right) =f\left( x\right) ,\text{ }%
x\in \left( 0,1\right)
\end{equation}

\begin{equation}
L_{k}u=\sum\limits_{i=0}^{m_{k}}\left[ \alpha _{ki}u^{\left( i\right)
}\left( 0\right) +\beta _{ki}u^{\left( i\right) }\left( 1\right)
+\sum\limits_{j=1}^{N_{k}}\delta _{kji}u^{\left( i\right) }\left(
x_{kj}\right) \right] =0,k=1,2,
\end{equation}%
where $A_{\lambda }=A+\lambda $, $A=A\left( x\right) $, $B=B\left( x\right) $
are linear operators in a Banach space $E,$ $a=a\left( x\right) $ is \ a
complex valued function, $\alpha _{ki},\beta _{ki},\delta _{kji}$ are
complex numbers, $x_{kj}\in \left( 0,1\right) $ and $\lambda $ is a spectral
parameter. Let we denote $\alpha _{km_{k}}$ and $\beta _{km_{k}}$by $\alpha
_{k}$ and $\beta _{k},$ respectively. Let\ $\omega _{1}=\omega _{1}\left(
x\right) ,\omega _{2}=\omega _{2}\left( x\right) $ be roots of the
characteristic equation $a\left( x\right) \omega ^{2}+1=0$ and $\ $%
\begin{equation*}
\eta =\eta \left( x\right) =\left\vert 
\begin{array}{cc}
\left( -\omega _{1}\right) ^{m_{1}}\alpha _{1} & \beta _{1}\omega
_{1}^{m_{1}} \\ 
\left( -\omega _{2}\right) ^{m_{2}}\alpha _{2} & \beta _{2}\omega
_{2}^{m_{2}}%
\end{array}%
\right\vert .
\end{equation*}

Function $u\in W_{p,\gamma }^{2}\left( 0,1;E\left( A\right) ,E\right) ,$ $%
L_{k}u=0$ satisfying the equation $\left( 7\right) $ a.e. on $\left(
0,1\right) $ is said to be solution of the problem $\left( 7\right) -\left(
8\right) .$

We say that the problem $\left( 7\right) -\left( 8\right) $ is $L_{p,\gamma
} $-separable, if for all $f\in L_{p,\gamma }\left( 0,1;E\right) $ there
exists a unique solution $u\in $ $W_{p,\gamma }^{2}\left( 0,1;E\left(
A\right) ,E\right) $ of the problem $\left( 7\right) -\left( 8\right) $ and
there exists a positive constant $C$ such that the coercive estimate holds 
\begin{equation*}
\left\Vert u^{\left( 2\right) }\right\Vert _{L_{p,\gamma }\left(
0,1;E\right) }+\left\Vert Au\right\Vert _{L_{p,\gamma }\left( 0,1;E\right)
}\leq C\left\Vert f\right\Vert _{L_{p,\gamma }\left( 0,1;E\right) }.
\end{equation*}

Let $Q$ denote the operator generated by BVP $\left( 7\right) -\left(
8\right) $ i.e.%
\begin{equation*}
D\left( Q\right) =W_{p,\gamma }^{2}\left( 0,1;E\left( A\right)
,E,L_{k}\right) ,\text{ }Qu=au^{\left( 2\right) }+Au+Bu^{^{\left( 1\right)
}}.
\end{equation*}

Let $I\left( E\left( A\right) ,E\right) $ denote the embedding operator from 
$E\left( A\right) $ to $E.$

\textbf{Condition 1. }Let the following conditions be satisfied:

(1) $E$ is an uniformly convex Banach space\ space with base and $\gamma \in
A_{p},$ $p\in \left( 1,\infty \right) $;

(2) $A$ is an $R$-positive in $E$ with $\varphi \in \left[ 0\right. \left.
\pi \right) ,$ $A\left( x\right) A^{-1}\left( \bar{x}\right) \in C\left( %
\left[ 0,1\right] ;B\left( E\right) \right) ,$ $\bar{x}\in \left( 0,1\right) 
$ and $BA^{\left( \frac{1}{2}-\mu \right) }\in L_{\infty }\left( 0,1;B\left(
E\right) \right) $ for $0<\mu <\frac{1}{2};$

(3) $-a\in S\left( \varphi _{1}\right) \cap \mathbb{C}\mathbf{/}\mathbb{R}%
_{-}$, $a\neq 0,$ $\eta \left( x\right) \neq 0,$ $\ 0\leq \varphi _{1}<\pi ,$
$\lambda \in S\left( \varphi _{2}\right) ,$ $\varphi _{1}+\varphi
_{2}<\varphi ,$;

Let $I=I\left( W_{p,\gamma }^{2}\left( 0,1;E\left( A\right) ,E\right)
,L_{p,\gamma }\left( 0,1;E\right) \right) $ denote the embedding operator%
\begin{equation*}
W_{p,\gamma }^{2}\left( 0,1;E\left( A\right) ,E\right) \rightarrow
L_{p,\gamma }\left( 0,1;E\right) .
\end{equation*}

In a similar way as in $\left[ \text{19, Theorem 3}\right] $ we obtain

\textbf{Theorem A}$_{6}.$Suppose the Condition1 holds. Then the problem $%
\left( 7\right) -\left( 8\right) $ for $f\in L_{p,\gamma }\left(
0,1;E\right) $, $\left\vert \arg \lambda \right\vert \leq \varphi $ and
sufficiently large $\left\vert \lambda \right\vert $ has a unique solution $%
u\in $ $W_{p,\gamma }^{2}\left( 0,1;E\left( A\right) ,E\right) $ and the
coercive uniform estimate holds

\begin{equation}
\sum\limits_{i=0}^{2}\left\vert \lambda \right\vert ^{1-\frac{i}{2}%
}\left\Vert u^{\left( i\right) }\right\Vert _{L_{p,\gamma }\left(
0,1;E\right) }+\left\Vert Au\right\Vert _{L_{p,\gamma }\left( 0,1;E\right)
}\leq M\left\Vert f\right\Vert _{L_{p,\gamma }\left( 0,1;E\right) }.
\end{equation}

Moreover from $\left[ 3\right] $ we have:

\textbf{Theorem A}$_{7}.$ Let $E$ be Banach spaces with base. Suppose the
operator $A$ is positive in $E$ and $A^{-1}\in \sigma _{\infty }\left(
E\right) .$ Assume that 
\begin{eqnarray*}
0 &\leq &\gamma <p-1,\text{ }1<p<\infty ,\text{ } \\
s_{j}\left( I\left( E\left( A\right) ,E\right) \right) &\sim &j^{-\frac{1}{%
\nu }},\text{ for some }\nu >0,\text{ }j=1,2,...,\infty .
\end{eqnarray*}%
Then the embedding $W_{p,\gamma }^{2}\left( 0,1;E\left( A\right) ,E\right)
\subset L_{p,\gamma }0,1;E$ is compact and 
\begin{equation*}
s_{j}\left( I\left( W_{p,\gamma }^{2}\left( 0,1;E\left( A\right) ,E\right)
,L_{p,\gamma }\left( 0,1;E\right) \right) \right) \sim j^{-\frac{2}{2\nu +1}%
}.
\end{equation*}

\textbf{Remake 3. }Really, Theorems A$_{6}$ and A$_{7}$ are proven under
condition that $E$ is an $\zeta $-convex Banach space. Since all uniformly
convex space is a $\zeta $-convex space i.e. is an UMD space, by applying $%
\left[ 3\right] $ we get the assertions.

By applying Theorem 3 and Theorems A$_{6},$ A$_{7}$ we obtain

\textbf{Theorem 4. }Suppose the Condition1 holds and 
\begin{equation*}
\text{ }s_{j}\left( I\left( E\left( A\right) ,E\right) \right) \sim j^{-%
\frac{1}{\nu }},\text{ for some }\nu >0,\text{ }j=1,2,...,\infty ;
\end{equation*}

Then:

(a) spectrum of the operator $Q$ is discrete;

(b) 
\begin{equation}
s_{j}\left( \left( Q+\lambda \right) ^{-1}\left( L_{p,\gamma }\left(
0,1;E\right) \right) \right) \sim j^{-\frac{2}{2\nu +1}}.
\end{equation}

(c) if $\varphi \leq \frac{\pi }{2q},$ $q>\nu +\frac{1}{2}$ then the system
of root functions of differential operator $Q$ is complete in $L_{p,\gamma
}\left( 0,1;E\right) .$

\textbf{Proof.} By virtue Theorem A$_{1},$ there exists a resolvent operator 
$\left( Q+\lambda \right) ^{-1}$ which is bounded from $L_{p,\gamma }\left(
0,1;E\right) $ to $W_{p,\gamma }^{2}\left( 0,1;E\left( A\right) ,E\right) .$
Moreover, by virtue of Theorem A$_{2}$ the embedding operator $I\left(
W_{p,\gamma }^{2}\left( 0,1;E\left( A\right) ,E\right) ,L_{p,\gamma }\left(
0,1;E\right) \right) $ is compact and 
\begin{equation}
s_{j}\left( I\left( W_{p,\gamma }^{2}\left( 0,1;E\left( A\right) ,E\right)
,L_{p,\gamma }\left( 0,1;E\right) \right) \right) \sim j^{-\frac{2}{2\nu +1}%
}.
\end{equation}

Since 
\begin{equation}
\left( Q+\lambda \right) ^{-1}\left( L_{p,\gamma }\left( 0,1;E\right)
\right) =\left( Q+\lambda \right) ^{-1}\left( L_{p,\gamma }\left(
0,1;E\right) ,W_{p,\gamma }^{2}\left( 0,1;E\left( A\right) ,E\right) \right)
\notag
\end{equation}%
\begin{equation*}
I\left( W_{p,\gamma }^{2}\left( 0,1;E\left( A\right) ,E\right) ,L_{p,\gamma
}\left( 0,1;E\right) \right)
\end{equation*}%
then from relations $\left( 11\right) $ and $\left( 12\right) $ we obtain
the assertions (a) and (b). Moreover, the estimate $\left( 9\right) $ and
the relation $\left( 11\right) $\ implies that operator $Q$ is positive in $%
L_{p,\gamma }\left( 0,1;E\right) $ and 
\begin{equation*}
\left( Q+\lambda \right) ^{-1}\in \tilde{\sigma}_{q}\left( L_{p,\gamma
}\left( 0,1;E\right) \right) ,\text{ for }q>\nu +\frac{1}{2}\text{ and }%
\lambda \in S\left( \varphi \right) .
\end{equation*}

By virtue of Remarke1, the above estimate implies 
\begin{equation}
\left( Q+\lambda \right) ^{-1}\in \sigma _{q}\left( L_{p,\gamma }\left(
0,1;E\right) \right) ,\text{ }q>\nu +\frac{1}{2}.
\end{equation}%
Then in view of the estimate $\left( 9\right) $, the relation $\left(
13\right) $ and by Theorem 3 we obtain the assertion (b).

Consider the following nonlocal BVP for degenerate DOE%
\begin{equation}
\left( L+\lambda \right) u=a\left( x\right) u^{\left[ 2\right] }\left(
x\right) +B\left( x\right) u^{^{\left[ 1\right] }}\left( x\right)
+A_{\lambda }\left( x\right) u\left( x\right) =f\left( x\right) ,\text{ }%
x\in \left( 0,1\right)
\end{equation}

\begin{equation*}
L_{k}u=\sum\limits_{i=0}^{m_{k}}\left[ \alpha _{ki}u^{\left[ i\right]
}\left( 0\right) +\beta _{ki}u^{\left[ i\right] }\left( 1\right)
+\sum\limits_{j=1}^{N_{k}}\delta _{kji}u^{\left[ i\right] }\left(
x_{kj}\right) \right] =0,k=1,2,
\end{equation*}%
where%
\begin{equation*}
u^{\left[ i\right] }=\left( x^{\gamma }\frac{d}{dx}\right) ^{i}u.
\end{equation*}
Let $O$ denote the operator generated by problem $\left( 14\right) $ and%
\begin{equation*}
W_{p}^{\left[ 2\right] }\left( 0,1;E_{0},E\right) =\left\{ u\in L_{p}\left(
0,1;E_{0}\right) \right. ,\text{ }u^{\left[ 2\right] }\in L_{p}\left(
0,1;E\right) ,
\end{equation*}

\begin{equation*}
\ \left\Vert u\right\Vert _{W_{p,\gamma }^{\left[ 2\right] }\left(
0,1;E_{0},E\right) }=\left\Vert u\right\Vert _{L_{p}\left( 0,1;E_{0}\right)
}+\left\Vert u^{\left[ 2\right] }\right\Vert _{L_{p}\left( 0,1;E\right)
}<\infty .\text{ }
\end{equation*}%
Theorem 4 implies the following result:

\textbf{Result 1. }Suppose all conditions of Theorem 4 are satisfies. Then
the assertions (a), (b) and (c) of Theorem 4 are hold for the operator $O$
in $L_{p}\left( 0,1;E\right) .$

Really, under the substitution%
\begin{equation*}
y=\int\limits_{0}^{x}z^{-\gamma }dz
\end{equation*}
the spaces $L_{p}\left( 0,1;E\right) ,$ $W_{p,\gamma }^{\left[ 2\right]
}\left( 0,1;E\left( A\right) ,E\right) $ are mapped isomorphically onto
spaces $L_{p,\tilde{\gamma}}\left( 0,b;E\right) $, $W_{p,\tilde{\gamma}%
}^{2}\left( 0,b;E\left( A\right) ,E\right) ,$ respectively, where $\tilde{%
\gamma}=\left[ \left( 1-\gamma \right) y\right] ^{\frac{1}{1-\gamma }}.$
Moreover, under this substitution the problem $\left( 14\right) $ is
transformed into a non degenerate problem $\left( 7\right) -\left( 8\right) $%
.

\begin{center}
\textbf{References}\ \ \ \ \ \ \ \ \ \ \ \ \ \ \ \ \ \ \ \ \ \ \ \ \ \ \ \ \
\ \ \ \ \ \ \ \ \ \ \ \ \ \ \ \ \ \ \ \ \ \ \ \ \ \ \ \ \ \ \ \ \ \ \ \ \ \
\ \ \ \ \ \ \ \ \ \ \ \ \ \ \ \ \ \ \ \ \ \ \ \ \ \ 
\end{center}

\begin{enumerate}
\item Amann H., Linear and quasi-linear equations,1, Birkhauser, 1995.

\item Agranovicn M. S., Spectral Boundary Value Problems in Lipschitz
Domains or Strongly Elliptic Systems in Banach Spaces $H_{p}^{\sigma }$ and $%
B_{p}^{\sigma }$, Functional Analysis and its Applications, 42 (4), (2008),
249-267.

\item Agarwal, R. P, Bohner, R., Shakhmurov, V. B, Maximal regular boundary
value problems in Banach-valued weighted spaces, Bound. Value Probl., 1
(2005), 9-42.

\item Burkholder D. L., A geometrical conditions that implies the existence
certain singular integral of Banach space-valued functions, Proc. conf.
Harmonic analysis in honor of Antonu Zigmund, Chicago, 1981,Wads Worth,
Belmont, (1983), 270-286.

\item Bourgain J., Some remarks on Banach spaces in which martingale
differences are unconditional, Arkiv Math. 21 (1983), 163-168.

\item Burgoyne J., Denseness of the generalized eigenvectors of a discrete
operator in a Banach space, J. Operator Theory, 33 (1995), 279--297.

\item Dore G. and Yakubov S., Semigroup estimates and non coercive boundary
value problems, Semigroup Form 60 (2000), 93-121.

\item Dunford N., Schwartz J. T., Linear operators. Parts 2$.$ Spectral
theory, Interscience, New York, 1963.

\item Denk R., Hieber M., Pr\"{u}ss J., $R$-boundedness, Fourier multipliers
and problems of elliptic and parabolic type, Mem. Amer. Math. Soc. 166
(2003), n.788.

\item Diestel, J., Geometry of Banach spaces-selected topics,
Springer-Verlag, Berlin-Heidelberg-New-York, 1975.

\item Gorbachuk V. I. and Gorbachuk M. L., Boundary value problems for
differential-operator equations, Naukova Dumka, Kiev, 1984.

\item K\"{o}nig, H., Eigenvalue Distribution of Compact Operators, Operator
Theory: Advances and Applications, vol. 16, Birkhauser, Basel etc., 1986.

\item Krein S. G., Linear differential equations in Banach space,
Providence, 1982.

\item Lunardi A., Analytic semigroups and optimal regularity in parabolic
problems, Birkhauser, 2003.

\item Markus A. S., Some criteria for the completeness of a system of root
vectors of a linear operator in a Banach space,\ Mat. Sb., 70 (112):4
(1966), 526--561.

\item Sobolevskii P. E., Inequalities coerciveness for abstract parabolic
equations, Dokl. Akad. Nauk. SSSR, (1964), 57(1), 27-40.

\item Shahmurov R., Solution of the Dirichlet and Neumann problems for a
modified Helmholtz equation in Besov spaces on an annuals, J. Differential
Equations, 249(3) (2010), 526-550.

\item Shakhmurov V. B., Imbedding theorems and their applications to
degenerate equations, Differential equations, 24 (4), (1988), 475-482.

\item Shakhmurov V. B., Coercive boundary value problems for regular
degenerate differential-operator equations, J. Math. Anal. Appl., 292 ( 2),
(2004), 605-620.

\item Shakhmurov V. B., Embedding theorems and\ maximal regular differential
operator equations in Banach-valued function spaces, J. Inequal. Appl.,
2(4)(2005), 329-345, 2( 4), (2005), 329-345.

\item Shakhmurov V. B., Separable anisotropic differential operators and
applications, J. Math. Anal. Appl. 2006, 327(2), 1182-1201.

\item Shakhmurov V. B., Embedding and maximal regular differential operators
in Banach-valued weighted spaces\textbf{, }Acta. Math. Sin., (Engl. Ser.),
(2012), 28 (9), 1883-1896.

\item Shahmurov R., On strong solutions of a Robin problem modeling heat
conduction in materials with corroded boundary, Nonlinear Anal. Real World
Appl., 13(1) (2011), 441-451.

\item Triebel H., Interpolation theory. Function spaces. Differential
operators.\ , North-Holland, Amsterdam, 1978.

\item Yakubov S., Completeness of root functions of regular differential
operators, Longman, Scientific and Technical, New York, 1994.

\item Yakubov S and Yakubov Ya., Differential-operator equations. Ordinary
and Partial \ Differential equations, Chapmen and Hall /CRC, Boca Raton,
2000.
\end{enumerate}

\bigskip

\end{document}